\documentstyle[11pt]{article}

\makeatletter

\@addtoreset{equation}{section}
\makeatother
\def\<{\langle}
\def\>{\rangle}

\newtheorem{lem}{Lemma}[section]
\newtheorem{theo}{Theorem}[section]
\newtheorem{rem}{Remark}[section]

\makeatletter
   
   \@addtoreset{equation}{section}
\makeatother

\begin{document}
\title{\bf Asymptotic Profile of Solutions to the Linearized Compressible Navier-Stokes Flow}
\author{{Ruy Coimbra Char\~ao}\\
{\small Department of Mathematics, Federal University of Santa
Catarina} \\ {\small 88040-270 Florianopolis, Brazil}\\Ryo Ikehata\thanks{Corresponding author: ikehatar@hiroshima-u.ac.jp} \\ {\small Department of Mathematics, Graduate School of Education, Hiroshima University} \\ {\small Higashi-Hiroshima 739-8524, Japan} }
\maketitle
\vspace{-0.8cm}
\begin{abstract}
We consider the asymptotic behavior as $t \to +\infty$ of the $L^{2}$-norm of the velocity of the linearized compressible Navier-Stokes equations in ${\bf R}^{n}$ ($n \geq 2$). As an application we shall study the optimality of the decay rate for the $L^{2}$-norm of the velocity by deriving the decay estimate from below as $t \to +\infty$. 
\end{abstract}

\section{Introduction}
\footnote[0]{Keywords and Phrases: Compressible Navier-Stokes equations; Cauchy problem; Asymptotic profiles; Weighted $L^{1}$-initial data; Low and high frequencies.}
\footnote[0]{2010 Mathematics Subject Classification. Primary 35Q30, 35B40; Secondary 76N99, 35C20.}
In this paper, we are concerned with the following linearized Compressible Navier-Stokes flow in ${\bf R}^{n}$ with $n \geq 2$:
\begin{equation}\label{eq1.1}
\hspace{-4.40cm}\rho_{t}(t,x) + \gamma{\rm div}v(t,x) = 0,\ \ \ (t,x)\in (0,\infty)\times {\bf R}^{n},
\end{equation}
\begin{equation}\label{eq1.2}
v_{t}(t,x)  - \alpha \Delta v(t,x) - \beta \nabla{\rm div}v(t,x) + \gamma \nabla \rho(t,x) = 0,\ \ \ (t,x)\in (0,\infty)\times {\bf R}^{n},
\end{equation}
\begin{equation}\label{eq1.3}
\hspace{-5.20cm}\rho(0,x)= \rho_{0}(x), \quad v(0,x) = v_{0}(x), \quad x\in {\bf R}^{n},
\end{equation}
where $\alpha $ and $\beta$, the viscosity  coefficient, are constants satisfying the thermodynamic restriction  $\alpha>0$ and  $\beta \geq 0$. The constant coefficient $\gamma$ is such that $\gamma > 0$  and 
$$v(t,x) :=\,\, ^{T}\Big(v_{1}(t,x), v_{2}(t,x), \cdots, v_{n}(t,x)\Big)$$ is the vector valued unknown  velocity of the fluid, $\rho(t,x)$ is a scalar valued unknown density of the fluid. Furthermore,  $v_{0}(x) :=\,\, ^{T}(v_{01}(x), v_{02}(x), \cdots, v_{0n}(x))$ and $\rho_{0}(x)$ are given initial data.


Concerning the $L^{p}$-$L^{q}$ estimates of the solution $(\rho(t,x), v(t,x))$ to the linearized NS equation (1.1)-(1.3), one should first mention a precise result due to Kobayashi-Shibata \cite{KS}, and in particular, they investigated the diffusion wave property of the solution in terms of the $L^{\infty}$-norm. The diffusion wave property was studied by Hoff-Zumbrum \cite{HZ,HZ-1} and Liu-Wang \cite{TW} to the (nonlinear)  compressible Navier-Stokes flow.  We should cite important closely related results  due to Brezina-Kagei \cite{BK}, Chowdhury-Ramaswamy \cite{CR}, Decklnik \cite{D} and the references therein as for the large time behavior of solutions to the compressible Navier-Stokes equations. In particular, we should mention the recent work by Ma-Wang \cite{MW}, which studied a large time asymptotic behavior of contact wave for the Cauchy problem of the one-dimensional compressible Navier-Stokes equations with restriction to zero viscosity. Furthermore, Matsumura-Nishida \cite{MN} studied the existence theory and time-asymptotic $L^{2}$ decay of solutions to the original (nonlinear) Navier-Stokes systems, and this theory is generalized to the more general hyperbolic-parabolic systems by Kawashima \cite{Kawa}. On the other hand, it is well-known that the density $\rho(t,x)$ becomes a solution of a corresponding viscoelastic equation, and therefore, the decay estimates due to Shibata \cite{shibata} and Ponce \cite{p} are very useful to investigate the problem (1.1)-(1.3). All these investigations are quite restricted to the decay property of the solutions in terms of $L^{p}$-norms. But, to the best of authors' knowledge, it seems that there still does not exist a research from the viewpoint that catches an asymptotic profile itself of the solutions to problem (1.1)-(1.3). Recently, in Liu-Noh \cite{LN} they have announced new interesting results about the asymptotic profile and its point-wise decay estimates of the Green function to the (linearized) Navier-Stokes System. Furthermore, quite recently Ikehata-Onodera \cite{IO} has caught the explicit profile of the "density" $\rho(t,x)$ as $t \to +\infty$ through the study of the viscoelastic equations, which is inspired from the previous results by \cite{IKM}. In this connection, the asymptotic profile of the viscoelastic equation was first discovered in the paper due to  Ikehata-Todorova-Yordanov  \cite{ITY}. The profile is so called the diffusion wave, which is popular in the field of the Navier-Stokes systems. The result of \cite{IO} reads as follows.

\begin{theo}\label{teo1.1}  Let $n \geq 2$. Then, it is true that there exist constants $C>0$ and $\eta > 0$ such that for $t > 0$ 
\begin{eqnarray*}
&&\int_{{\bf R}^{n}}\Big| \;\hat{\rho}(t,\xi) - \Big[-(i\xi)\cdot P_{0}e^{-\frac{(\alpha + \beta)t\vert\xi\vert^{2}}{2}}\frac{\sin(\gamma t\vert\xi\vert)}{\vert\xi\vert} + Q_{0}e^{-\frac{(\alpha + \beta)t\vert\xi\vert^{2}}{2}}\cos(\gamma t\vert\xi\vert)\;\Big]\;\Big|^{2}\,d\xi\nonumber\\
&\leq& Ct^{-\frac{n}{2}-1}\Vert \rho_{0}\Vert_{1,1}^{2}+Ct^{-\frac{n}{2}-1}\vert P_{0}\vert^{2} + Ct^{-\frac{n}{2}-1}\vert Q_{0}\vert^{2} + Ct^{-\frac{n}{2}-1}(\;\sum_{j=1}^{n}\Vert v_{0j}\Vert_{1,1}^{2}\;)\nonumber\\
&+& Ce^{-\eta t}\;(\;\Vert {\rm div}v_{0}\Vert^{2} + \Vert \rho_{0}\Vert^{2} + \vert Q_{0}\vert^{2} + \vert P_{0}\vert^{2}\;),
\end{eqnarray*}
for $t \gg 1$, where
\[P_{0} := (P_{01}, P_{02},\cdots, P_{0n}), \quad P_{0j} := \int_{{\bf R}^{n}}v_{0j}(x)dx\,\, (j = 1,2,\cdots,n), \quad Q_{0} := \int_{{\bf R}^{n}}\rho_{0}(x)dx\]
and $\hat{\rho}(t,\xi)$ is the usual Fourier transform of  $\rho(t,x)$.
\end{theo}
However,  until now we still did not know about the explicit profile of the velocity $v(t,x)$ of the fluid . 

The main purpose of this note is to announce the exact profile of the velocity $v(t,x)$ as $t \to +\infty$ (see Lemma 3.1 below).
\begin{theo} \label{teo1.2}
Let $n \geq 2$. Then, it is true that there exists generous constants $C > 0$ and $\eta > 0$ such that
\begin{eqnarray*}
&&\int_{{\bf R}^{n}} \Big| \hat{v}(t,\xi) - P_{0}e^{-\alpha\vert\xi\vert^{2}t} + \frac{\xi(\xi\cdot P_{0})}{\vert\xi\vert^{2}}e^{-\alpha\vert\xi\vert^{2}t}\\
&+&(i\xi)e^{-\frac{(\alpha + \beta)\vert\xi\vert^{2}t}{2}} \frac{\sin(\gamma t\vert\xi\vert)}{\vert\xi\vert}Q_{0} - \frac{\xi(\xi\cdot P_{0})}{\vert\xi\vert^{2}}e^{-\frac{(\alpha + \beta)\vert\xi\vert^{2}t}{2}}\cos(\gamma t\vert\xi\vert)\Big|^{2}d\xi\\
&\leq& C\;(\;\vert P_{0}\vert^{2} + \vert Q_{0}\vert^{2} + \sum_{j=1}^{n}\Vert v_{0j}\Vert_{1,1}^{2} + \Vert\rho_{0}\Vert_{1,1}^{2})t^{-\frac{n}{2}-1} + Ce^{-\eta t}(\;\Vert v_{0}\Vert^{2} + \Vert \rho_{0}\Vert^{2}\;),
\end{eqnarray*}
for large $t \gg 1$, where $C > 0$ depends only on $\gamma$, $\alpha$, $\beta$, and $\hat{v}(t,\xi)$ is the Fourier transform of $v(t,x)$.
\end{theo}

As an application, one has the optimal decay estimate of the $L^{2}$-norm of the velocity $v(t,x)$ for the compressible fluid. 
\begin{theo}\label{teo1.3}  
Let $n \geq 2$. Then, it is true that there exist constants $C_{j} >0$ {\rm (}$j = 1,2${\rm )} such that for $t \gg 1$ 
\[C_{1}\Big(\;\vert P_{0}\vert + \vert Q_{0}\vert\;\Big)t^{-\frac{n}{4}} \leq \Vert v(t,\cdot)\Vert \leq C_{2}\Big(\;\vert P_{0}\vert + \vert Q_{0}\vert + \sum_{j=1}^{n}\Vert v_{0j}\Vert_{1,1} + \Vert\rho_{0}\Vert_{1,1} + \Vert v_{0}\Vert + \Vert \rho_{0}\Vert\;\Big)t^{-\frac{n}{4}},\]
provided that $\vert Q_{0}\vert \ne 0$ and $\vert P_{0}\vert/\vert Q_{0}\vert \ll 1$.
\end{theo}
\begin{rem}\,{\rm  It is still open to show the optimality above in the case when the assumptions $\vert Q_{0}\vert \ne 0$ and $\vert P_{0}\vert/\vert Q_{0}\vert \ll 1$ do not hold. However, the assumption $\vert P_{0}\vert/\vert Q_{0}\vert \ll 1$ means that
$\vert P_{0}\vert / \vert Q_{0}\vert < \displaystyle{\frac{1}{1+C_n}}$ with $C_n>0$ depending on the dimension $n$ and the coefficients $\alpha, \; \beta$, which can be calculated explicitly.} 

\vspace{0.5cm}

{\footnotesize {\bf Notation.} Throughout this paper, $\| \cdot\|_q$ stands for the usual $L^q({\bf R}^{n})$-norm. For simplicity of notations, in particular, we use $\| \cdot\|$ instead of $\| \cdot\|_2$. 
Furthermore, we set
\[f \in L^{1,\gamma}({\bf R}^{n}) \Leftrightarrow f \in L^{1}({\bf R}^{n}),\, \Vert f\Vert_{1,\gamma} := \int_{{\bf R}^{n}}(1+\vert x\vert^{\gamma})\vert f(x)\vert dx < +\infty, \quad \gamma > 0,\]
\[\Vert f\Vert := \sqrt{\Vert f_{1}\Vert^{2} + \cdots \Vert f_{n}\Vert^{2}}\]
for $f = (f_{1}, \cdots, f_{n}) \in (L^{2}({\bf R}^{n}))^{n}$.\\

On the other hand, we denote the Fourier transform $\hat{\phi}(\xi)$ of the function $\phi(x)$ by

\[{\cal F}(\phi)(\xi) := {\cal F}_{x \to \xi}(\phi)(\xi) := \hat{\phi}(\xi) := \frac{1}{(2\pi)^{n/2}}\int_{{\bf R}^{n}}e^{-ix\cdot\xi}\phi(x)dx,\]
and we denote by ${\cal F}_{x \to \xi}^{-1}$ its usual inverse Fourier transform, where $i := \sqrt{-1}$, and $x\cdot\xi = \displaystyle{\sum_{j=1}^{n}}x_{j}\xi_{j}$ for $x = (x_{1},\cdots,x_{n})$ and $\xi = (\xi_{1},\cdots,\xi_{n})$. We also use the notation
\[v_{t}=\frac{\partial u}{\partial t}, \quad \Delta = \sum^n_{j=1}\frac{\partial^2}{\partial x_j^2}, \quad \nabla f = \nabla f(x) := (\frac{\partial f}{\partial x_{1}}, \cdots,\frac{\partial f}{\partial x_{n}}),\]\\
and the notation $div v$ means the usual divergence of the vector valued function $v$.\\
}

\end{rem}
\begin{rem}\,{\rm One of our advantage is the ease of the method of our proof as compared with that of \cite{KS} when we catch the exact profile of the solution {\rm (}$\rho(t,x), v(t,x)${\rm )}. We use the method from \cite{Ik-3} or \cite{Ik-0} when one deals with the low frequency part of the Fourier transformed velocity, while in the high frequency region we shall rely on the method due to \cite{LIC} combined with the Haraux-Komornik inequality (see Lemma 2.4 below), which is a version of the energy method in the Fourier space. The latter part seems much different from known techniques about the compressible Navier-Stokes equations.}
\end{rem}



\begin{rem}\,{\rm 
In the case when (for example) $n = 3$, it follows from the study due to \cite[p. 395]{LN} that
\[{\cal F}_{x\to \xi}^{-1}\left(P_{0j}\frac{\xi_{j}\xi_{k}}{\vert\xi\vert^{2}}(\cos(\gamma t\vert\xi\vert)-1)e^{-\frac{(\alpha+\beta)\vert\xi\vert^{2}t}{2}}\right)(x)\]
\[ = P_{0j}\frac{\partial^{2}}{\partial x_{j}\partial x_{k}}\frac{\gamma^{2}}{4\pi}\int_{0}^{t}\int_{\vert y\vert = 1}G(\frac{\alpha+\beta}{2}t,x + \gamma\tau y)dS_{y}d\tau,\]
and
\[{\cal F}_{{\bf x}\to {\bf \xi}}^{-1}\left(P_{0j}\frac{\xi_{j}\xi_{k}}{\vert\xi\vert^{2}}(e^{-\frac{(\alpha+\beta)\vert\xi\vert^{2}t}{2}}-e^{-\alpha\vert\xi\vert^{2}t})\right)(x)\]
\[ = P_{0j}\frac{\partial^{2}}{\partial x_{j}\partial x_{k}}\frac{\alpha-\beta}{2}\int_{0}^{t}\int_{\vert y\vert = 1}G(\min\{\frac{\alpha+\beta}{2},\alpha\}t + \frac{\vert\alpha-\beta\vert}{2}\tau, x + \gamma\tau y)dS_{ y}d\tau,\]
where $j,k = 1,2,\cdots,3$, and
\[G(t,x) = (4\pi t)^{-\frac{3}{2}}e^{-\frac{\vert x\vert^{2}}{4t}}\]
is the $3$-dimensional Gauss kernel. Therefore, if one can know the decomposition below:
\[P_{0}e^{-\alpha\vert\xi\vert^{2}t} - \frac{\xi(\xi\cdot P_{0})}{\vert\xi\vert^{2}}e^{-\alpha\vert\xi\vert^{2}t} - (i\xi)e^{-\frac{(\alpha + \beta)\vert\xi\vert^{2}t}{2}}\frac{\sin(\gamma t\vert\xi\vert)}{\vert\xi\vert}Q_{0} + \frac{\xi(\xi\cdot P_{0})}{\vert\xi\vert^{2}}e^{-\frac{(\alpha + \beta)\vert\xi\vert^{2}t}{2}}\cos(\gamma t\vert\xi\vert)\]
\[= P_{0}e^{-\alpha\vert\xi\vert^{2}t} -(i\xi)e^{-\frac{(\alpha + \beta)\vert\xi\vert^{2}t}{2}}\frac{\sin(\gamma t\vert\xi\vert)}{\vert\xi\vert}Q_{0} + \frac{\xi(\xi\cdot P_{0})}{\vert\xi\vert^{2}}e^{-\frac{(\alpha + \beta)\vert\xi\vert^{2}t}{2}}\cos(\gamma t\vert\xi\vert) - \frac{\xi(\xi\cdot P_{0})}{\vert\xi\vert^{2}}e^{-\alpha\vert\xi\vert^{2}t}\]

\begin{equation}
= P_{0}e^{-\alpha\vert\xi\vert^{2}t} -(i\xi)e^{-\frac{(\alpha + \beta)\vert\xi\vert^{2}t}{2}}\frac{\sin(\gamma t\vert\xi\vert)}{\vert\xi\vert}Q_{0}
\end{equation}
\[+ \frac{\xi(\xi\cdot P_{0})}{\vert\xi\vert^{2}}\left(e^{-\frac{(\alpha + \beta)\vert\xi\vert^{2}t}{2}}-e^{-\alpha\vert\xi\vert^{2}t}\right) + \frac{\xi(\xi\cdot P_{0})}{\vert\xi\vert^{2}}\left(\cos(\gamma t\vert\xi\vert)-1\right)e^{-\frac{(\alpha + \beta)\vert\xi\vert^{2}t}{2}}
\]
one can get a precise profile of the velocity $v(t,x)$. This is left to the readers' check. Note that as for the one of profiles for (1.4) part in asymptotic sense, one can refer the reader to \cite[p. 2167]{Ik-0} in order to obtain the formula: for each $k = 1,2,3$ (cf. \cite[Lemma 3.3]{KS})
\[{\cal F}_{x \to \xi}^{-1}\left((i\xi_{k})\hat{w}(t,\cdot)\hat{h}(t,\cdot)\right)(x) = \frac{\partial}{\partial x_{k}}(w(t,\cdot)*h(t,\cdot))(x)\]
\[= C_{0}(\alpha + \beta)^{-\frac{3}{2}}t^{-\frac{1}{2}}\frac{\partial}{\partial x_{k}}\int_{\vert z\vert = 1}e^{-\frac{\vert x + tz\vert^{2}}{2(\alpha + \beta)t}}dS_{z}\]
with some constant $C_{0} > 0$. Here the function $w(t,x)$ is the fundamental solution to the free wave equation
\[
w_{tt}(t,x)  - \gamma^{2}\Delta w(t,x) = 0,\ \ \ (t,x)\in (0,\infty)\times {\bf R}^{3},
\]
\[
w(0,x)= 0, \quad w_{t}(0,x) = \delta(x), \quad x\in {\bf R}^{3},
\]
where $\delta(x)$ is the usual Dirac measure, and 
\[h(t,x) := (\alpha+\beta)^{-\frac{3}{2}}t^{-\frac{3}{2}}e^{-\frac{\vert x\vert^{2}}{2(\alpha + \beta)t}}.\]
It is easy to check
\[\hat{w}(t,\xi) = \frac{\sin(\gamma t\vert\xi\vert)}{\vert\xi\vert},\quad \hat{h}(t,\xi) = e^{-\frac{(\alpha+\beta)t\vert\xi\vert^{2}}{2}}.\]
}
\end{rem}

In the rest of this paper we shall give a proof of Theorem \ref{teo1.2} in Section 2, and in Section 3 we will prove Theorem \ref{teo1.3}.


\section{Proof of Theorem \ref{teo1.2}.}

\hspace{0.6cm}In this section, we shall prove Theorem \ref{teo1.2} based on a device due to \cite{IO}.
We first apply the Fourier transform to problem (\ref{eq1.1})--(\ref{eq1.3}). Then the problem (\ref{eq1.1})--(\ref{eq1.3}) can be reduced to the following system of ordinary differential equations with frequency parameter $\xi = (\xi_{1},\xi_{2},\cdots,\xi_{n}) \in {\bf R}_{\xi}^{n}$
\begin{equation}\label{eq2.1}
\hspace{-4.6cm}\hat{\rho}_{t}(t,\xi) + i\gamma\xi\cdot\hat{v}(t,\xi) = 0,\ \ \ (t,\xi)\in (0,\infty)\times {\bf R}_{\xi}^{n},
\end{equation}
\begin{equation}\label{eq2.2}
\hat{v}_{t}(t,\xi)  + \alpha\vert\xi\vert^{2}\hat{v}(t,\xi) + \beta\xi(\xi\cdot\hat{v}(t,\xi)) + i\gamma\xi\hat{\rho}(t,\xi) = 0,\ \ \ (t,\xi)\in (0,\infty)\times {\bf R}_{\xi}^{n},
\end{equation}
\begin{equation}\label{eq2.3}
\hspace{-5.5cm}\hat{\rho}(0,\xi)= \hat{\rho}_{0}(\xi), \quad \hat{v}(0,\xi) = \hat{v}_{0}(\xi), \quad \xi \in {\bf R}_{\xi}^{n}.
\end{equation}

We first prove the following lemma in the low frequency region in the case when $0 < \vert\xi\vert \ll 1$.

\begin{lem}\label{lem2.1} 
Let $n \geq 2$. Then, there exists a small $\delta_{0} > 0$ such that
\begin{eqnarray*}
&&\int_{\vert\xi\vert \leq \frac{\delta_{0}}{\sqrt{2}}}\Big\vert\;\hat{v}(t,\xi) - P_{0}e^{-\alpha\vert\xi\vert^{2}t} + \frac{\xi(\xi\cdot P_{0})}{\vert\xi\vert^{2}}e^{-\alpha\vert\xi\vert^{2}t}\\
&+&(i\xi)e^{-(\alpha + \beta)\vert\xi\vert^{2}t/2} \frac{\sin(\gamma t\vert\xi\vert)}{\vert\xi\vert}Q_{0} - \frac{\xi(\xi\cdot P_{0})}{\vert\xi\vert^{2}}e^{-(\alpha + \beta)\vert\xi\vert^{2}t/2}\cos(\gamma t\vert\xi\vert)\;\Big\vert^{2}d\xi\\
&\leq & C\;\Big(\;\vert P_{0}\vert^{2} + \vert Q_{0}\vert^{2} + \sum_{j=1}^{n}\Vert v_{0j}\Vert_{1,1}^{2} + \Vert\rho_{0}\Vert_{1,1}^{2}\;\Big)\;t^{-\frac{n}{2}-1},
\end{eqnarray*}
for large $t \gg 1$, where $C > 0$ is a generous constant depending only on $\gamma$, $\alpha$, $\beta$, and so on.
\end{lem}

In order to prove Lemma \ref{lem2.1} above, let us solve (\ref{eq2.1})-(\ref{eq2.3}) directly from the viewpoint of the velocity $v(t,x)$ under the condition that $0 < \vert\xi\vert \leq \displaystyle{\frac{\delta_{0}}{\sqrt{2}}}$ with small $\delta_{0} > 0$ by basing on the result due to \cite[(2.9)]{KS}. In this case we get

\begin{eqnarray}\label{eq2.4}
\hat{v}(t,\xi) &=& e^{-\alpha\vert\xi\vert^{2}t}\hat{v}_{0}(\xi) - (i\gamma\xi)\Big(\;\frac{e^{\sigma_{1}t}-e^{\sigma_{2}t}}{\sigma_{1}-\sigma_{2}}\;\Big)\hat{\rho}_{0}(\xi)
\nonumber \\
&+& \Big(\;\frac{\sigma_{1}e^{\sigma_{1}t}-\sigma_{2}e^{\sigma_{2}t}}{\sigma_{1}-\sigma_{2}}- e^{-\alpha\vert\xi\vert^{2}t}\;\Big)\frac{\xi(\xi\cdot\hat{v}_{0}(\xi))}{\vert\xi\vert^{2}},
\end{eqnarray}
provided that $0 < \vert\xi\vert \leq \displaystyle{\frac{\delta_{0}}{\sqrt{2}}}$, where $\sigma_{j} \in {C}$ ($j = 1,2$) are given by the following expressions
\[\sigma_{1}= { \sigma_{1}(\xi)} = \frac{-b\vert\xi\vert^{2}+i\vert\xi\vert\sqrt{4a-b^{2}\vert\xi\vert^{2}}}{2}, \quad \sigma_{2}={ \sigma_{2}(\xi)} = \frac{-b\vert\xi\vert^{2}-i\vert\xi\vert\sqrt{4a-b^{2}\vert\xi\vert^{2}}}{2},\]
with
\[a := \gamma^{2}, \quad b := (\alpha + \beta), \quad \delta_{0} := \frac{2\sqrt{a}}{b} = \frac{2\gamma}{\alpha + \beta}.\]

Now let us use an idea  which was  introduced in \cite{Ik-3}. We consider the following  decomposition of the initial data
\begin{equation}\label{eq2.5}
\hat{v}_{0j}(\xi) = A_{0j}(\xi) -iB_{0j}(\xi) + P_{0j}, \quad (j = 1,2,\cdots, n),
\end{equation}
\begin{equation}\label{eq2.6}
\hat{\rho}_{0}(\xi) = A_{\rho}(\xi) -iB_{\rho}(\xi) + Q_{0}, 
\end{equation}
where
\[A_{\rho}(\xi) := \int_{{\bf R}^{n}}\Big(\;\cos(x\cdot\xi)-1\;\Big)\rho_{0}(x) dx, \quad B_{\rho}(\xi) := \int_{{\bf R}^{n}}\sin(x\cdot\xi)\rho_{0}(x) dx,\]
\[A_{0j}(\xi) := \int_{{\bf R}^{n}}\Big(\;\cos(x\cdot\xi)-1\;\Big)v_{0j}(x) dx, \quad B_{0j}(\xi) := \int_{{\bf R}^{n}}\sin(x\cdot\xi)v_{0j}(x) dx, \quad(j = 1,2,\cdots, n),\]
\[Q_{0} := \int_{{\bf R}^{n}}\rho_{0}(x)dx,\quad P_{0j} :=  \int_{{\bf R}^{n}}v_{0j}(x)dx, \quad(j =1,2,\cdots, n),\]
where $v_{0j}$ are the components of  the initial velocity  $v_0$.

Since we can write 
\[(i\gamma)\xi\cdot\hat{v}_{0}(\xi) = (i\gamma)\sum_{j=1}^{n}\xi_{j}\cdot(\;A_{0j}-iB_{0j} + P_{0j}\;) =: (i\gamma\xi)\cdot\Big(A_{0}(\xi) -iB_{0}(\xi) + P_{0}\Big),\]
where
\[A_{0}(\xi) -iB_{0}(\xi) + P_{0}\]
\[:= \Big(A_{01}(\xi), A_{02}(\xi),\cdots,A_{0n}(\xi)\Big) - i\Big(B_{01}(\xi), B_{02}(\xi),\cdots,B_{0n}(\xi)\Big) + \Big(P_{01},P_{02},\cdots,P_{0n}\Big),\]
one has the following expression for the velocity $\hat{v}(t,\xi)$
\begin{eqnarray}\label{eq2.7}
\hat{v}(t,\xi) &=& e^{-\alpha\vert\xi\vert^{2}t}\Big[A_{0}(\xi) -iB_{0}(\xi) + P_{0}\Big] - (i\gamma\xi)(\frac{e^{\sigma_{1}t}-e^{\sigma_{2}t}}{\sigma_{1}-\sigma_{2}})\Big[A_{\rho}(\xi) -iB_{\rho}(\xi) + Q_{0}\Big] 
\nonumber\\
&+& \Big(\;\frac{\sigma_{1}e^{\sigma_{1}t}-\sigma_{2}e^{\sigma_{2}t}}{\sigma_{1}-\sigma_{2}}-e^{-\alpha\vert\xi\vert^{2}t}\;\Big)\frac{\xi\Big[\;\xi\cdot\Big(A_{0}(\xi) -iB_{0}(\xi) + P_{0}\Big)\;\Big]}{\vert\xi\vert^{2}},
\end{eqnarray}
for all $\xi$ satisfying $0 < \vert\xi\vert \leq \displaystyle{\frac{\delta_{0}}{\sqrt{2}}}$.

Now, it is easy to check that
\begin{equation}\label{eq2.8}
\frac{e^{\sigma_{1}t}-e^{\sigma_{2}t}}{\sigma_{1}-\sigma_{2}} = 2\frac{e^{-bt\vert\xi\vert^{2}/2}\sin(\frac{t\vert\xi\vert\sqrt{4a-b^{2}\vert\xi\vert^{2}}}{2})}{\vert\xi\vert\sqrt{4a-b^{2}\vert\xi\vert^{2}}},
\end{equation}
and
\begin{equation}\label{eq2.9}
\frac{\sigma_{1}e^{\sigma_{1}t}-\sigma_{2}e^{\sigma_{2}t}}{\sigma_{1}-\sigma_{2}} = -\frac{b\vert\xi\vert e^{-bt\vert\xi\vert^{2}/2}\sin(\frac{t\vert\xi\vert\sqrt{4a-b^{2}\vert\xi\vert^{2}}}{2})}{\sqrt{4a-b^{2}\vert\xi\vert^{2}}} + e^{-bt\vert\xi\vert^{2}/2}\cos(\frac{t\vert\xi\vert\sqrt{4a-b^{2}\vert\xi\vert^{2}}}{2}).
\end{equation}
\noindent
So, it follows from (\ref{eq2.7}), (\ref{eq2.8}) and (\ref{eq2.9}) that
\begin{eqnarray}\label{eq2.10}
\hat{v}(t,\xi)&=&P_{0}e^{-\alpha\vert\xi\vert^{2}t} - \xi(\xi\cdot P_{0})\frac{b e^{-bt\vert\xi\vert^{2}/2}\sin(\frac{t\vert\xi\vert\sqrt{4a-b^{2}\vert\xi\vert^{2}}}{2})}{\vert\xi\vert\sqrt{4a-b^{2}\vert\xi\vert^{2}}}\nonumber\\
&+& \frac{\xi(\xi\cdot P_{0})}{\vert\xi\vert^{2}}e^{-bt\vert\xi\vert^{2}/2}\cos(\frac{t\vert\xi\vert\sqrt{4a-b^{2}\vert\xi\vert^{2}}}{2})
\\
&-&e^{-\alpha\vert\xi\vert^{2}t}\frac{\xi(\xi\cdot P_{0})}{\vert\xi\vert^{2}} -2(i\gamma\xi)\frac{e^{-b\vert\xi\vert^{2}t/2}\sin(\frac{t\vert\xi\vert\sqrt{4a-b^{2}\vert\xi\vert^{2}}}{2})}{\vert\xi\vert\sqrt{4a-b^{2}\vert\xi\vert^{2}}}Q_{0} + E_{0}(t,\xi),\nonumber
\end{eqnarray}
where
\begin{eqnarray}\label{eq2.11}
E_{0}(t,\xi)& :=&  e^{-\alpha\vert\xi\vert^{2}t}\Big[\;A_{0}(\xi) -iB_{0}(\xi)\;\Big] - (i\gamma\xi)(\;\frac{e^{\sigma_{1}t}-e^{\sigma_{2}t}}{\sigma_{1}-\sigma_{2}}\;)\Big[\;A_{\rho}(\xi) -iB_{\rho}(\xi)\;\Big] \nonumber\\
&+& \Big(\;\frac{\sigma_{1}e^{\sigma_{2}t}-\sigma_{2}e^{\sigma_{1}t}}{\sigma_{1}-\sigma_{2}}-e^{-\alpha\vert\xi\vert^{2}t}\;\Big)\frac{\xi\Big(\;\xi\cdot\Big[\;A_{0}(\xi) -iB_{0}(\xi)\;\Big]\;\Big)}{\vert\xi\vert^{2}}.
\end{eqnarray}

 Now, applying the mean value theorem for $|\xi| \leq \displaystyle{\frac{\delta_{0}}{\sqrt{2}}}$ it follows that
\begin{equation}\label{eq2.12}
2\frac{\sin(\frac{t\vert\xi\vert\sqrt{4a-b^{2}\vert\xi\vert^{2}}}{2})}{\vert\xi\vert\sqrt{4a-b^{2}\vert\xi\vert^{2}}} = \frac{2}{\sqrt{4a-b^{2}\vert\xi\vert^{2}}}\frac{\sin(\sqrt{a}t\vert\xi\vert)}{\vert\xi\vert} + t(\;\frac{\sqrt{4a-b^{2}\vert\xi\vert^{2}}-2\sqrt{a}}{\sqrt{4a-b^{2}\vert\xi\vert^{2}}}\;)\cos(\varepsilon(t,\xi))
\end{equation}
and
\begin{equation}\label{eq2.13}
\cos(\frac{t\vert\xi\vert\sqrt{4a-b^{2}\vert\xi\vert^{2}}}{2}) = \cos(\sqrt{a}t\vert\xi\vert) - t\vert\xi\vert(\;\frac{\sqrt{4a-b^{2}\vert\xi\vert^{2}}-2\sqrt{a}}{2}\;)\sin(\eta(t,\xi)),
\end{equation}
where
\[\varepsilon(t,\xi) := \frac{t\vert\xi\vert\sqrt{4a-b^{2}\vert\xi\vert^{2}}}{2}\theta + \sqrt{a}t\vert\xi\vert(1-\theta'),\]
\[\eta(t,\xi) := \frac{t\vert\xi\vert\sqrt{4a-b^{2}\vert\xi\vert^{2}}}{2}\theta' + \sqrt{a}t\vert\xi\vert(1-\theta''),\]
for some $\theta'$ and  $\theta'' \in (0,1)$.
\noindent
Furthermore, using again the mean value theorem it follows that 
\begin{equation}\label{eq2.14}
\frac{2}{\sqrt{4a-b^{2}\vert\xi\vert^{2}}} = \frac{1}{\sqrt{a}} + \frac{2b^{2}\theta\vert\xi\vert^{2}}{(4a-b^{2}\theta^{2}\vert\xi\vert^{2})\sqrt{4a-b^{2}\theta^{2}\vert\xi\vert^{2}}}, \quad \theta \in (0,1).
\end{equation}
\noindent
 Then,  from (\ref{eq2.10})--(\ref{eq2.14}) in the case when $0 < \vert\xi\vert \leq \displaystyle{\frac{\delta_{0}}{\sqrt{2}}}$ we find the following useful expression for the Fourier transform of the velocity $v(t,x)$ 
 
\begin{eqnarray}\label{eq2.15}
\hat{v}(t,\xi)&=&P_{0}e^{-\alpha\vert\xi\vert^{2}t} - \frac{\xi(\xi\cdot P_{0})}{\vert\xi\vert^{2}}e^{-\alpha\vert\xi\vert^{2}t} - (i\xi)e^{-b\vert\xi\vert^{2}t/2} \frac{\sin(\gamma t\vert\xi\vert)}{\vert\xi\vert}Q_{0} \nonumber\\
&-& \frac{b}{2}\xi(\xi\cdot P_{0})e^{-b\vert\xi\vert^{2}t/2}\frac{\sin(\gamma t\vert\xi\vert)}{\gamma\vert\xi\vert}  +  \frac{\xi(\xi\cdot P_{0})}{\vert\xi\vert^{2}}e^{-b\vert\xi\vert^{2}t/2}\cos(\gamma t\vert\xi\vert)\nonumber \\
&+& \,E_{0}(t,\xi) -  \frac{\xi(\xi\cdot P_{0})}{\vert\xi\vert^{2}}e^{-b\vert\xi\vert^{2}t/2}(t\vert\xi\vert)\frac{\sqrt{D}-2\gamma}{2}\sin(\eta(t,\xi))\\
&-&(i\gamma\xi)Q_{0}e^{-b\vert\xi\vert^{2}t/2}\sin(\gamma t\vert\xi\vert)\frac{2b^{2}\theta\vert\xi\vert}{\sqrt{D_{\theta}^{3}}} - (i\gamma\xi)Q_{0}t e^{-b\vert\xi\vert^{2}t/2}(\frac{\sqrt{D}-2\gamma}{\sqrt{D}})\cos(\varepsilon(t,\xi))\nonumber\\
&-& b^{3}\xi(\xi\cdot P_{0})\frac{\theta\vert\xi\vert}{\sqrt{D_{\theta}^{3}}}e^{-b\vert\xi\vert^{2}t/2}\sin(\gamma t\vert\xi\vert) - b\xi(\xi\cdot P_{0})(\frac{t}{2})e^{-b\vert\xi\vert^{2}t/2}(\frac{\sqrt{D}-2\gamma}{\sqrt{D}})\cos(\varepsilon(t,\xi)), \nonumber
\end{eqnarray}
where $D := 4a -b^{2}\vert\xi\vert^{2}$ and $D_{\theta} := 4a-b^{2}\theta^{2}\vert\xi\vert^{2}$. 

In order to estimate the remainder term we set
\begin{eqnarray*}
E_{1}(t,\xi) &:=&  -\frac{\xi(\xi\cdot P_{0})}{\vert\xi\vert^{2}}e^{-b\vert\xi\vert^{2}t/2}(t\vert\xi\vert)\frac{\sqrt{D}-2\gamma}{2}\sin(\eta(t,\xi)),\\
E_{2}(t,\xi)  &:=& b^{3}\xi(\xi\cdot P_{0})\frac{\theta\vert\xi\vert}{\sqrt{D_{\theta}^{3}}}e^{-b\vert\xi\vert^{2}t/2}\sin(\gamma t\vert\xi\vert), \\
E_{3}(t,\xi)  &:=&  -b\xi(\xi\cdot P_{0})\frac{t}{2}e^{-b\vert\xi\vert^{2}t/2}(\frac{\sqrt{D}-2\gamma}{\sqrt{D}})\cos(\varepsilon(t,\xi)),\\
E_{4}(t,\xi)  &:=& -(i\gamma\xi)Q_{0}e^{-b\vert\xi\vert^{2}t/2}\sin(\gamma t\vert\xi\vert)\frac{2b^{2}\theta\vert\xi\vert}{\sqrt{D_{\theta}^{3}}},\\
E_{5}(t,\xi)  &:=& - (i\gamma\xi)Q_{0}t e^{-b\vert\xi\vert^{2}t/2}(\frac{\sqrt{D}-2\gamma}{\sqrt{D}})\cos(\varepsilon(t,\xi)),\\
E_{6}(t,\xi)  &:=& - \frac{b}{2}\xi(\xi\cdot P_{0})e^{-b\vert\xi\vert^{2}t/2}\frac{\sin(\gamma t\vert\xi\vert)}{\gamma\vert\xi\vert}.\\
\end{eqnarray*}

In fact,  for the profile that we consider,  $\displaystyle{\sum_{j=0}^{6}}E_{j}(t,\xi)$ is a remainder term. Moreover, the next calculations will show that it is the best  choice. This implies that  the Fourier transform of the velocity $v(t,x)$ is given by  
\begin{eqnarray}\label{eq2.16}
\hat{v}(t,\xi)&=& P_{0}e^{-\alpha\vert\xi\vert^{2}t} - \frac{\xi(\xi\cdot P_{0})}{\vert\xi\vert^{2}}e^{-\alpha\vert\xi\vert^{2}t} - (i\xi)e^{-b\vert\xi\vert^{2}t/2} \frac{\sin(\gamma t\vert\xi\vert)}{\vert\xi\vert}Q_{0}\nonumber\\
&+&\frac{\xi(\xi\cdot P_{0})}{\vert\xi\vert^{2}}e^{-b\vert\xi\vert^{2}t/2}\cos(\gamma t\vert\xi\vert) + \displaystyle{\sum_{j=0}^{6}}E_{j}(t,\xi)
\end{eqnarray}
in the low frequency zone $0 < \vert\xi\vert \leq \displaystyle{\frac{\delta_{0}}{\sqrt{2}}}$.

Now, let us estimate all quantities $E_{j}(t,\xi)$ in terms of $L^{2}({\bf R}_{\xi}^{n})$-norm in order to make sure the fact that $\displaystyle{\sum_{j=0}^{6}}{ E}_{j}(t,\xi)$ is the remainder term.  For this ends, we first prepare the following elementary  helpful estimate
\begin{equation}\label{eq2.17}
\vert \sqrt{4a-b^{2}\vert\xi\vert^{2}}-2\sqrt{a}\vert = \vert\frac{b^{2}\vert\xi\vert^{2}}{\vert \sqrt{4a-b^{2}\vert\xi\vert^{2}}+2\sqrt{a}}\vert \leq \frac{b^{2}\vert\xi\vert^{2}}{2\sqrt{a}},
\end{equation}
{which hold for  $0 < \vert\xi\vert \leq \displaystyle{\frac{\delta_{0}}{\sqrt{2}}} < \delta_0=\displaystyle{\frac{2\sqrt{a}}{b}}$.} { We also note that  
$$4a-b^2\theta^2|\xi|^2 \geq  4a-b^2|\xi|^2 \geq 2a $$
for $0 < \vert\xi\vert \leq \displaystyle{\frac{\delta_{0}}{\sqrt{2}}}$  and  $0< \theta< 1$.

Based on the inequalities (\ref{eq2.16}) and (\ref{eq2.17}), by using the Schwarz inequality we can proceed all estimates below except for  $E_{0}(t,\xi)$.}

\begin{eqnarray}\label{eq2.18}
\int_{\vert\xi\vert \leq \delta_{0}/\sqrt{2}}\vert E_{1}(t,\xi)\vert^{2}d\xi &\leq& t^{2}\vert P_{0}\vert^{2}\int_{\vert\xi\vert \leq \delta_{0}/\sqrt{2}}\vert\xi\vert^{2}e^{-b\vert\xi\vert^{2}t}\frac{b^{4}\vert\xi\vert^{4}}{4a}d\xi\\
&\leq& \frac{b^{4}}{4a}\vert P_{0}\vert^{2}t^{2}\int_{\vert\xi\vert \leq \delta_{0}/\sqrt{2}}\vert\xi\vert^{6}e^{-b\vert\xi\vert^{2}t}d\xi \leq \frac{b^{4}}{4a}\vert P_{0}\vert^{2}t^{-\frac{n}{2}-1},\nonumber
\end{eqnarray}

\begin{eqnarray}\label{eq2.19}
\int_{\vert\xi\vert \leq \delta_{0}/\sqrt{2}}\vert E_{2}(t,\xi)\vert^{2}d\xi &\leq& b^{6}\vert P_{0}\vert^{2}\int_{\vert\xi\vert\leq\delta_{0}/\sqrt{2}}\vert\xi\vert^{6}e^{-b\vert\xi\vert^{2}t}\frac{1}{\vert 4a-b^{2}\theta^{2}\vert\xi\vert^{2} \vert^{3}}d\xi\\
&\leq& \frac{b^{6}\vert P_{0}\vert^{2}}{({\ 2a})^{3}}\int_{\vert\xi\vert \leq \delta_{0}/\sqrt{2}}\vert\xi\vert^{6}e^{-b\vert\xi\vert^{2}t}d\xi \leq \frac{b^{6}\vert P_{0}\vert^{2}}{({ 2a})^{3}}t^{-\frac{n}{2}-3},\nonumber
\end{eqnarray}

\begin{eqnarray}\label{eq2.20}
\int_{\vert\xi\vert \leq \delta_{0}/\sqrt{2}}\vert E_{3}(t,\xi)\vert^{2}d\xi 
&\leq & b^{2}t^{2}\vert P_{0}\vert^{2}\int_{\vert\xi\vert \leq \delta_{0}/\sqrt{2}}\vert\xi\vert^{4}e^{-b\vert\xi\vert^{2}t}\frac{\vert \sqrt{4a-b^{2}\vert\xi\vert^{2}}-2\sqrt{a}\vert^{2}}{4a-b^{2}\vert\xi\vert^{2}}d\xi\nonumber\\
&\leq & b^{2}\vert P_{0}\vert^{2}t^{2}\int_{\vert\xi\vert \leq \delta_{0}/\sqrt{2}}\vert\xi\vert^{4}e^{-b\vert\xi\vert^{2}t}\frac{b^{4}\vert\xi\vert^{4}}{4a(4a-b^{2}\vert\xi\vert^{2})}d\xi \\
&\leq & \frac{b^{6}\vert P_{0}\vert^{2}t^{2}}{{ 8a}^{2}}\int_{\vert\xi\vert \leq \delta_{0}/\sqrt{2}}\vert\xi\vert^{8}e^{-b\vert\xi\vert^{2}t}d\xi \leq \frac{b^{6}\vert P_{0}\vert^{2}t^{2}}{{ 8a}^{2}}t^{-\frac{n}{2}-2},\nonumber
\end{eqnarray}

\begin{eqnarray}\label{eq2.21}
\int_{\vert\xi\vert \leq \delta_{0}/\sqrt{2}}\vert E_{4}(t,\xi)\vert^{2}d\xi
& \leq & 4\gamma^{2}b^{4}\vert Q_{0}\vert^{2}\int_{\vert\xi\vert \leq \delta_{0}/\sqrt{2}}\vert\xi\vert^{4}\frac{e^{-b\vert\xi\vert^{2}t}}{(4a-b^{2}\theta^{2}\vert\xi\vert^{2})^{3}}d\xi\nonumber
\\
&\leq& \frac{4\gamma^{2}b^{4}\vert Q_{0}\vert^{2}}{\vert 3a\vert^{3}}\int_{\vert\xi\vert \leq \delta_{0}/\sqrt{2}}\vert\xi\vert^{4}e^{-b\vert\xi\vert^{2}t}d\xi\\ 
&\leq& \frac{4\gamma^{2}b^{4}\vert Q_{0}\vert^{2}}{\vert { 2a}\vert^{3}}t^{-\frac{n}{2}-2},
\nonumber
\end{eqnarray}

\begin{eqnarray}\label{eq2.22}
\int_{\vert\xi\vert \leq \delta_{0}/\sqrt{2}}\vert E_{5}(t,\xi)\vert^{2}d\xi 
&\leq& \gamma^{2}t^{2}\vert Q_{0}\vert^{2}\int_{\vert\xi\vert \leq \delta_{0}/\sqrt{2}}\vert\xi\vert^{2}e^{-b\vert\xi\vert^{2}t}\frac{\vert \sqrt{4a-b^{2}\vert\xi\vert^{2}}-2\sqrt{a}\vert^{2}}{4a-b^{2}\vert\xi\vert^{2}}d\xi\nonumber\\
&\leq& \frac{\gamma^{2}t^{2}\vert Q_{0}\vert^{2}}{{ 2a}}\int_{\vert\xi\vert \leq \delta_{0}/\sqrt{2}}\vert\xi\vert^{2}e^{-b\vert\xi\vert^{2}t}\frac{b^{4}\vert\xi\vert^{4}}{4a}d\xi\\
&\leq& \frac{b^{4}\gamma^{2}\vert Q_{0}\vert^{2}}{{\ 8a}^{2}}t^{2}\int_{\vert\xi\vert \leq \delta_{0}/\sqrt{2}}\vert\xi\vert^{6}e^{-b\vert\xi\vert^{2}t}d\xi \leq \frac{b^{4}\gamma^{2}\vert Q_{0}\vert^{2}}{{ 8a}^{2}}t^{-\frac{n}{2}-1},\nonumber
\end{eqnarray}

\begin{equation}\label{eq2.23}
\int_{\vert\xi\vert \leq \delta_{0}/\sqrt{2}}\vert E_{6}(t,\xi)\vert^{2}d\xi \leq \frac{b^{2}}{4\gamma^{2}}\vert P_{0}\vert^{2}\int_{\vert\xi\vert \leq \delta_{0}/\sqrt{2}}\vert\xi\vert^{2}e^{-b\vert\xi\vert^{2}t}d\xi \leq \frac{b^{2}}{4\gamma^{2}}\vert P_{0}\vert^{2}t^{-\frac{n}{2}-1}.
\end{equation}

In the above estimates we also  have used the elementary estimate 
\begin{eqnarray}\label{elem-est}
 \int_{|\xi|\leq \frac{\delta_0}{\sqrt{2}}} |\xi|^k e^{-\lambda|\xi|^2t} \;d\xi \leq C_n\;t^{-\frac{n+k}{2}} \;, \qquad t>0,
 \end{eqnarray}
for $\lambda >0$, $\;k+n>0$ where $C_n$ is a positive constant depending on $\delta_0$, $\lambda$, $k$ and the dimension $n$.

In order to estimate $E_{0}(t,\xi)$, we prepare the following simple lemma which plays an essential role in this note. This idea has its origin in \cite[Lemma 3.1]{Ik-3}.

\begin{lem}\label{lem2.2}
 Let $n \geq 1$. Then it holds that
\[\vert A_{\rho}(\xi)\vert \leq L\vert\xi\vert\Vert \rho_{0}\Vert_{1,1},\quad \vert A_{0j}(\xi)\vert \leq L\vert\xi\vert\Vert v_{0j}\Vert_{1,1},\quad (j = 1,2,\cdots,n),\]
\[\vert B_{\rho}(\xi)\vert \leq M\vert\xi\vert\Vert \rho_{0}\Vert_{1,1},\quad \vert B_{0j}(\xi)\vert \leq M\vert\xi\vert\Vert v_{0j}\Vert_{1,1},\quad (j = 1,2,\cdots,n),\]
  for all $\xi \in {\bf R}^{n}$, where
\[L := \sup_{\theta \ne 0}\frac{\vert 1-\cos\theta\vert}{\vert\theta\vert} < +\infty, \quad M := \sup_{\theta \ne 0}\frac{\vert \sin\theta\vert}{\vert\theta\vert} < +\infty.\]
We note that $L<1$ and $M=1$.
\end{lem}
\noindent

The estimate for $E_{0}(t,\xi)$ based on Lemma \ref{lem2.2} is crucial in this paper.\\ 

In fact, it follows  from (\ref{eq2.8}), (\ref{eq2.9}), (\ref{eq2.11}), Lemma \ref{lem2.2} and the Schwarz inequality that
\begin{eqnarray*}
\int_{\vert\xi\vert \leq \delta_{0}/\sqrt{2}}&&\hspace{-0.8cm}\vert E_{0}(t,\xi)\vert^{2}d\xi 
\;\;\leq \;\;C\int_{\vert\xi\vert \leq \delta_{0}/\sqrt{2}}e^{-2\alpha\vert\xi\vert^{2}t}\Big[\;\vert A_{0}(\xi)\vert^{2} + \vert B_{0}(\xi)\vert^{2}\;\Big]d\xi\\
&+& C\gamma^{2}\int_{\vert\xi\vert \leq \delta_{0}/\sqrt{2}}e^{-b\vert\xi\vert^{2}t}\vert\xi\vert^{2}\Big[\;\vert A_{\rho}(\xi)\vert^{2} + \vert B_{\rho}(\xi)\vert^{2}\;\Big]\frac{1}{\vert\xi\vert^{2}(4a-b^{2}\vert\xi\vert^{2})}d\xi\\
&+&C\int_{\vert\xi\vert \leq \delta_{0}/\sqrt{2}}\Big(\frac{b^{2}\vert\xi\vert^{2}e^{-b\vert\xi\vert^{2}t}}{4a-b^{2}\vert\xi\vert^{2}}+ e^{-b\vert\xi\vert^{2}t} +e^{-2\alpha\vert\xi\vert^{2}t} \Big)\frac{\vert\xi\vert^{4}\Big[\;\vert A_{0}(\xi)\vert^{2} + \vert B_{0}(\xi)\vert^{2}\;\Big]}{\vert\xi\vert^{4}}d\xi\\
&\leq&  C(L^{2}+ M^{2})(\;\sum_{j=1}^{n}\Vert v_{0j}\Vert_{1,1}^{2}\;)\int_{\vert\xi\vert \leq \delta_{0}/\sqrt{2}}\vert\xi\vert^{2}e^{-2\alpha\vert\xi\vert^{2}t}d\xi \\
&+& \frac{C\gamma^{2}(L^{2}+ M^{2})}{ 2a}\Vert\rho_{0}\Vert_{1,1}^{2}\int_{\vert\xi\vert \leq \delta_{0}/\sqrt{2}}\vert\xi\vert^{2}e^{-b\vert\xi\vert^{2}t}d\xi \\
&+& \frac{Cb^{2}}{2a}\int_{\vert\xi\vert \leq \delta_{0}/\sqrt{2}}\vert\xi\vert^{2}e^{-b\vert\xi\vert^{2}t}\Big[\;\vert A_{0}(\xi)\vert^{2} + \vert B_{0}(\xi)\vert^{2}\;\Big]d\xi\\
&+& C\int_{\vert\xi\vert \leq \delta_{0}/\sqrt{2}}e^{-\min\{b,2\alpha\}\vert\xi\vert^{2}t}\Big[\;\vert A_{0}(\xi)\vert^{2} + \vert B_{0}(\xi)\vert^{2}\;\Big]d\xi,
\end{eqnarray*}
 because of the fact that $4a-b^2|\xi|^2 \geq 2a$ for $|\xi|\leq \frac{\delta_0}{\sqrt{2}}$.

The previous inequality  and (\ref{elem-est}) imply that
\begin{eqnarray}\label{eq2.24}
\int_{\vert\xi\vert \leq \delta_{0}/\sqrt{2}}&&\hspace{-0.8cm}\vert E_{0}(t,\xi)\vert^{2}d\xi 
\;\leq \; C(L^{2}+ M^{2})(\;\sum_{j=1}^{n}\Vert v_{0j}\Vert_{1,1}^{2}\;)t^{-\frac{n}{2}-1} + \frac{C\gamma^{2}(L^{2}+ M^{2})}{2a}\Vert\rho_{0}\Vert_{1,1}^{2}t^{-\frac{n}{2}-1}
\nonumber\\
&+&\frac{Cb^{2}(L^{2}+ M^{2})}{2a}(\;\sum_{j=1}^{n}\Vert v_{0j}\Vert_{1,1}^{2}\;)\int_{\vert\xi\vert \leq \delta_{0}/\sqrt{2}}\vert\xi\vert^{4}e^{-b\vert\xi\vert^{2}t}d\xi \nonumber\\
&+&C(L^{2}+ M^{2})(\;\sum_{j=1}^{n}\Vert v_{0j}\Vert_{1,1}^{2}\;)\int_{\vert\xi\vert \leq \delta_{0}/\sqrt{2}}\vert\xi\vert^{2}e^{-\min\{b,2\alpha\}\vert\xi\vert^{2}t}d\xi\\
&\leq&  C(L^{2}+ M^{2})(\;\sum_{j=1}^{n}\Vert v_{0j}\Vert_{1,1}^{2}\;)t^{-\frac{n}{2}-1} + \,\frac{C\gamma^{2}(L^{2}+ M^{2})}{2a}\Vert\rho_{0}\Vert_{1,1}^{2}t^{-\frac{n}{2}-1}
\nonumber\\
&+&\frac{Cb^{2}(L^{2}+ M^{2})}{2a}(\;\sum_{j=1}^{n}\Vert v_{0j}\Vert_{1,1}^{2}\;)t^{-\frac{n}{2}-2} + \,C(L^{2}+ M^{2})(\;\sum_{j=1}^{n}\Vert v_{0j}\Vert_{1,1}^{2}\;)t^{-\frac{n}{2}-1},
\nonumber 
\end{eqnarray}
where the constants after the last inequality  depend on the dimension $n$, and they can be calculated explicitly.

 The expression for $\hat{v}(t,\xi)$ in (\ref{eq2.16}) combined with the estimates   (\ref{eq2.18})--(\ref{eq2.23}) and (\ref{eq2.24}) above and the fact that $b=\alpha+\beta$, imply the statement of Lemma \ref{lem2.1} in the low frequency region. \par\hfill$\Box$
\par

\vspace{0.3cm}
Finally in this section, we shall derive decay estimates in the high frequency region by relying on a special multiplier method in the Fourier space introduced in   Char\~ao-daLuz-Ikehata \cite{RCR2} (see also\cite{LIC}).

\begin{lem}\label{lem2.3} Let $n \geq 2$. Then, there exists a constant $\eta > 0$ such that
\[\int_{\vert\xi\vert \geq \frac{\delta_{0}}{\sqrt{2}}}\vert\hat{v}(t,\xi)\vert^{2}d\xi \leq C(\Vert v_{0}\Vert^{2} + \Vert \rho_{0}\Vert^{2})e^{-\eta t} \quad (t \gg 1),\]
where $\delta_{0} := \displaystyle{\frac{2\gamma}{\alpha + \beta}} = \displaystyle{\frac{2\sqrt{a}}{b}} > 0$.
\end{lem}
{\it Proof.} We prove the lemma for any arbitrary high frequency zone, that is, for any fixed number $\lambda_0>0$ in place of $\displaystyle{\frac{\delta_0}{\sqrt{2}}}$ by using the multiplier method combined with a simple version of the Komornik lemma. This result is stronger than that of Kobayashi--Shibata \cite{KS}, which proves similar results in high frequency zone only with $\lambda_0 >>1$.

In the proof of this lemma, in order to simplify the notation, we use 
$\hat{\rho}$ and  $ \hat{v}$ in place of $\hat{\rho}(t,\xi)$ and  $ \hat{v}(t,\xi)$ and
the same for $\hat{v}_t$ and $\hat{\rho}_t$, respectively. Moreover, we can use $\hat{v}(S)$ and 
$\hat{\rho}(S)$ in place of  $\hat{v}(S,\xi)$ and  $\hat{\rho}(S,\xi)$, respectively.

Multiply equation (\ref{eq2.1}) by $\hat{\rho}$ and equation (\ref{eq2.2}) by $\hat{v}$. Then we obtain 
\begin{eqnarray*}
\frac{d}{dt}\Big( \frac{|\hat{\rho}|^2+|\hat{v}|^2}{2}\Big) + \alpha \vert\xi\vert^{2}|\hat{v}|^2+ \beta |\xi \cdot \hat{v}|^2 + 2i Re(\gamma \bar{\hat{\rho}} \hat{v} \cdot \xi)=0 .
\end{eqnarray*}
\noindent
The above identity says that $Re(\gamma \bar{\hat{\rho}} \hat{v} \cdot \xi)=0 $. Thus, integrating the identity above on $[S, T]$ we get 
\begin{eqnarray}\label{eq3.00}
\Big( \frac{|\hat{\rho}|^2+|\hat{v}|^2}{2}\Big)_S^T + \alpha \int_S^T \vert\xi\vert^{2}|\hat{v}|^2dt+ \beta \int_S^T  |\xi \cdot \hat{v}|^2 dt=0, 
\end{eqnarray}
for all $ 0<S<T$ and  $ \xi \in {\bf R}^{n}$. Then it follows that
\begin{eqnarray}\label{eq3.0}
\alpha \int_S^T \vert\xi\vert^{2}|\hat{v}|^2dt \leq  \frac{|\hat{\rho}(S)|^2+|\hat{v}(S)|^2}{2}, \qquad 0<S<T, \qquad \xi \in {\bf R}^{n}. 
\end{eqnarray}
\noindent
In particular, 
$$E(T)+ \alpha \int_S^T \vert\xi\vert^{2}|\hat{v}|^2dt + \beta \int_S^T  |\xi \cdot \hat{v}|^2 dt = E(S), \qquad 0<S<T, \;\xi \in {\bf R}^{n}, $$
where $$E(t)=E(t,\xi) = \frac{|\hat{\rho}(t, \xi)|^2+|\hat{v}(t, \xi)|^2}{2}.$$
\noindent
Multiplying the equation (\ref{eq2.2}) by $\xi \bar{\hat{\rho}}$ we obtain 
\begin{eqnarray*}
\xi \cdot \Big( \hat{v}_t\bar{\hat{\rho}} + \alpha \vert\xi\vert^{2} \hat{v}\bar{\hat{\rho}} + \beta \xi (\xi \cdot \hat{v})\bar{\hat{\rho}} + i \gamma \xi |\hat{\rho}|^2 \Big) =0 
\end{eqnarray*}
or 
\begin{eqnarray*}
(\xi \cdot \hat{v}_t)\bar{\hat{\rho}} + \alpha \vert\xi\vert^{2} (\xi \cdot \hat{v})\bar{\hat{\rho}} + \beta |\xi|^2 (\xi \cdot \hat{v})\bar{\hat{\rho}} + i \gamma |\xi|^2 |\hat{\rho}|^2=0 .
\end{eqnarray*}

Now, multiplying  by $i=\sqrt{-1}$ it results
\begin{eqnarray*}
i(\xi \cdot \hat{v}_t)\bar{\hat{\rho}} + i\alpha \vert\xi\vert^{2} (\xi \cdot \hat{v})\bar{\hat{\rho}} + i\beta |\xi|^2 (\xi \cdot \hat{v})\bar{\hat{\rho}} = \gamma |\xi|^2 |\hat{\rho}|^2.
\end{eqnarray*}
\noindent
Then, by integrating it over $[S,T]$ one gets 
\begin{eqnarray}\label{eq3.1}
\gamma \int_S^T |\xi|^2 |\hat{\rho}|^2 dt = \int_S^T \Big[ i(\xi \cdot \hat{v}_t)\bar{\hat{\rho}} + i\alpha \vert\xi\vert^{2} (\xi \cdot \hat{v})\bar{\hat{\rho}} + i\beta |\xi|^2 (\xi \cdot \hat{v})\bar{\hat{\rho}}\Big] dt,
\end{eqnarray}
for $0<S<T$ and $\xi \in {\bf R}^{n}$.
\noindent
Since $\gamma >0$, the equation (\ref{eq2.1}) says that  $i(\xi \cdot \hat{v}) = - \displaystyle{\frac{1}{\gamma}}\hat{\rho}_t$. So, by substituting this fact in identity (\ref{eq3.1}) it follows that
\begin{eqnarray*}
\gamma \int_S^T |\xi|^2 |\hat{\rho}|^2 dt & = \displaystyle{\int_S^T} \Big[ i(\xi \cdot \hat{v}_t)\bar{\hat{\rho}} -\alpha |\xi|^2 \frac{1}{\gamma} \hat{\rho}_t\bar{\hat{\rho}}-\beta |\xi|^2 \frac{1}{\gamma} \hat{\rho}_t \bar{\hat{\rho}}\Big]dt\\
& = \displaystyle{\int_S^T} \Big[ i(\xi \cdot \hat{v}_t)\bar{\hat{\rho}} -\Big(\frac{\alpha}{\gamma} + \frac{\beta}{\gamma} \Big)|\xi|^2 \frac{d}{dt}\Big( |\hat{\rho}|^2\Big) \Big]dt\\
& = \displaystyle{\int_S^T}  i(\xi \cdot \hat{v}_t)\bar{\hat{\rho}} dt -\Big[ \frac{\alpha+ \beta}{\gamma}|\xi|^2   |\hat{\rho}(t)|^2\Big]_S^T.
\end{eqnarray*}
\noindent
Thus, one has
\begin{eqnarray} \label{eq3.2}
\gamma \int_S^T |\xi|^2 |\hat{\rho}|^2 dt  \leq \int_S^T  i(\xi \cdot \hat{v}_t)\bar{\hat{\rho}} dt + \frac{\alpha+ \beta}{\gamma}|\xi|^2 |\hat{\rho}(S)|^2,
\end{eqnarray}
for $0<S<T$ and $\xi \in {\bf R}^{n}$.

In order to estimate the integral in the right hand side of (\ref{eq3.2}) we can use two possibilities:  integration by parts or taking the conjugate of the equation (\ref{eq2.2}) and multiplying by $\hat{v}_t$. If we employ  the second option we get 
\begin{eqnarray*}
|\hat{v}_t|^2+\alpha |\xi|^2\bar{\hat{v}} \cdot \hat{v}_t + \beta (\xi \cdot \hat{v}_t)(\xi \cdot \bar{\hat{v}}) + i \gamma \xi \bar{\hat{\rho}} \cdot \hat{v}_t=0
\end{eqnarray*} 
or
\begin{eqnarray*}
|\hat{v}_t|^2+\alpha |\xi|^2 \frac{d}{dt}\frac{|\hat{v}|^2}{2} + \beta \frac{d}{dt}\frac{|\xi \cdot \hat{v}|^2}{2} + i \gamma (\xi \cdot \hat{v}_t) \bar{\hat{\rho}} =0.
\end{eqnarray*} 
\noindent
By integrating the above identity on $[S,T]$ we have 
\begin{eqnarray*}
\int_S^T |\hat{v}_t|^2\, dt +\Big[ \alpha |\xi|^2 \frac{|\hat{v}|^2}{2} + \beta \frac{|\xi \cdot \hat{v}|^2}{2}\Big]_S^T + \int_S^T i \gamma (\xi \cdot \hat{v}_t) \bar{\hat{\rho}}\, dt=0,
\end{eqnarray*} 
which implies 
\begin{eqnarray*}
\int_S^T i \gamma (\xi \cdot \hat{v}_t) \bar{\hat{\rho}}\, dt \leq 
\Big[ \alpha |\xi|^2 \frac{|\hat{v}|^2}{2} + \beta \frac{|\xi \cdot \hat{v}|^2}{2}\Big]_{t=S} \leq C_{\alpha, \beta} |\xi|^2 |\hat{v}(S)|^2, 
\end{eqnarray*}
where $\; C_{\alpha, \beta}>0$ is a constant depending only on $\alpha$ or $\beta$. Then it follows that
\begin{eqnarray} \label{eq3.3}
\int_S^T i (\xi \cdot \hat{v}_t) \bar{\hat{\rho}}\, dt \leq \frac{ C_{\alpha, \beta}}{ \gamma} |\xi|^2 |\hat{v}(S)|^2, \qquad  0<S<T, \qquad \xi \in {\bf R}^{n}.
\end{eqnarray}
\noindent
By substituting (\ref{eq3.3}) into (\ref{eq3.2}) one has obtained
\begin{eqnarray} \label{eq3.4}
\gamma \int_S^T |\xi|^2 |\hat{\rho}|^2 dt  \leq  \frac{C_{\alpha, \beta}}{\gamma} |\xi|^2 |\hat{v}(S)|^2 + \frac{\alpha+ \beta}{\gamma}|\xi|^2 |\hat{\rho}(S)|^2,
\end{eqnarray}
for $\xi \in {\bf R}^{n}$ and $0<S<T$. By combining (\ref{eq3.4}) and (\ref{eq3.0}) we arrive at 
\begin{eqnarray*} 
\int_S^T |\xi|^2 \Big[ |\hat{v}|^2+ |\hat{\rho}|^2\Big] dt 
 &\leq&  \frac{C_{\alpha, \beta}}{\gamma^2} |\xi|^2 |\hat{v}(S)|^2  
+\frac{\alpha+ \beta}{\gamma^2}|\xi|^2 |\hat{\rho}(S)|^2\\
 &+& \frac{|\xi|^2}{\lambda_0^2}\frac{|\hat{\rho}(S)|}{2\alpha}+\frac{|\xi|^2}{\lambda_0^2} \frac{|\hat{v}(S)|}{2\alpha},
\end{eqnarray*}
for all $|\xi| \geq \lambda_0>0$ and $0<S<T$. Then, we have obtained the following important estimate 
\begin{eqnarray}\label{eq3.5} 
\int_S^T  \Big[ |\hat{v}|^2+ |\hat{\rho}|^2\Big] dt  \leq  C_{\alpha, \beta, \gamma, \lambda_0} \Big[ |\hat{v}(S)|^2 + |\hat{\rho}(S)|^2\Big],
\end{eqnarray}
for all $|\xi| \geq \lambda_0>0$ and $0<S<T$, where $C_{\alpha, \beta, \gamma, \lambda_0} >0$ is a positive constant depending on $\alpha, \; \beta, \; \gamma \;\; \mbox{and}\;\;\delta_0$.

By using the definition of the energy for the system (\ref{eq2.1})--(\ref{eq2.2}) in the Fourier space
$$E(t,\xi): =|\hat{v}(t, \xi)|^2 + |\hat{\rho}(t, \xi)|^2,$$
(\ref{eq3.5}) implies that
\begin{eqnarray}\label{eq3.6}
\int_S^T \int_{|\xi|\leq \lambda_0} E(t,\xi)d\xi dt \leq C_{\alpha, \beta, \gamma, \lambda_0}\int_{|\xi|\leq \lambda_0} E(S,\xi)d\xi,
\end{eqnarray}
for $0<S<T<\infty$. 

Now, if we define the energy in high frequency zone in the Fourier space by
$$E_h(t):=\int_{|\xi|\leq \lambda_0} E(t,\xi)d\xi,$$
the estimate (\ref{eq3.6}) says that 
\begin{eqnarray}\label{eq3.7}
\int_S^{\infty} E_h(t) dt \leq C_{\alpha, \beta, \gamma, \lambda_0}E_h(S),
\end{eqnarray} 
for all $S \geq 0$.\par\hfill$\Box$

To get the final estimate for the energy just defined above on the high frequency zone, we use a simple version of the following well-known Haraux--Komornik lemma.

\begin{lem}\label{lem3.2}
Let $\,E :[0,+\infty) \rightarrow [0,+\infty)\,$ be a non-increasing
function and assume that there exists  a constant 
$\,T_0 > 0$ such that
$$\displaystyle{ \int_S^\infty} E(t)\,dt \leq T_0 \; E(S),$$
for all $S\; \geq 0$. Then, it is true that
$$\,\,E(t) \leq
E(0)\,e^{1-\frac{t}{T_0}}$$
for all $\,t \geq T_0$.
\end{lem} 

{\it Proof of Lemma 2.3 completed.}
In order to finalize the proof of Lemma \ref{lem2.3}, we note the energy in the high frequency region $E_h(t)$ is a non-increasing function of $t$ due to the identity (\ref{eq3.00}). Then, we can combine the estimate (\ref{eq3.7}) and Lemma \ref{lem3.2} to conclude that 
$$E_h(t) \leq CE_{h}(0)e^{-\eta t}  {\leq C \Big(\; ||v_0||^2 + ||\rho_0||^2 \;\Big) e^{-\eta t}}$$
for  $\eta= \displaystyle{\frac{1}{C_{\alpha, \beta, \gamma, \lambda_0}}}$ and $t \geq T_0:= C_{\alpha, \beta, \gamma, \lambda_0}$, where $C$ is a positive constant depending only on the coefficients of the system (\ref{eq1.1})--(\ref{eq1.2}) and $\lambda_0$. In particular, the above inequality proves the desired lemma with $\lambda_0=\displaystyle{\frac{\delta_0}{\sqrt{2}}}$.

\par\hfill
$\Box$


{\it Proof of Theorem 1.2.} The proof of Theorem 1.2 is a direct consequence of Lemmas 2.1 and 2.3.\par\hfill $\Box$\par

\section{Proof of Theorem 1.3.}

In this section, we shall give a proof of Theorem 1.3. For this ends it suffices to get the following lemma because of the Plancherel theorem together with Theorem 1.2 and the useful inequality:
\begin{eqnarray}
\Vert v(t,\cdot)\Vert&=& \Vert\hat{v}(t,\cdot)\Vert
\geq \Big\Vert P_{0}e^{-\alpha\vert\xi\vert^{2}t} - \frac{\xi(\xi\cdot P_{0})}{\vert\xi\vert^{2}}e^{-\alpha\vert\xi\vert^{2}t}\nonumber\\
&-& (i\xi)e^{-(\alpha + \beta)\vert\xi\vert^{2}t/2} \frac{\sin(\gamma t\vert\xi\vert)}{\vert\xi\vert}Q_{0} + \frac{\xi(\xi\cdot P_{0})}{\vert\xi\vert^{2}}e^{-(\alpha + \beta)\vert\xi\vert^{2}t/2}\cos(\gamma t\vert\xi\vert)\Big\Vert\nonumber\\
&-&\Big\Vert \hat{v}(t,\xi) - P_{0}e^{-\alpha\vert\xi\vert^{2}t} + \frac{\xi(\xi\cdot P_{0})}{\vert\xi\vert^{2}}e^{-\alpha\vert\xi\vert^{2}t}\nonumber\\
&+&\,\,(i\xi)e^{-(\alpha + \beta)\vert\xi\vert^{2}t/2} \frac{\sin(\gamma t\vert\xi\vert)}{\vert\xi\vert}Q_{0} - \frac{\xi(\xi\cdot P_{0})}{\vert\xi\vert^{2}}e^{-(\alpha + \beta)\vert\xi\vert^{2}t/2}\cos(\gamma t\vert\xi\vert)\Big\Vert\nonumber\\
&=& \Big\Vert (i\xi)\frac{\sin(\gamma t\vert\xi\vert)}{\vert\xi\vert}Q_{0}e^{-(\alpha + \beta)\vert\xi\vert^{2}t/2} - \{P_{0} - \frac{\xi(\xi\cdot P_{0})}{\vert\xi\vert^{2}}\}e^{-\alpha\vert\xi\vert^{2}t}\nonumber\\ 
&-&\frac{\xi(\xi\cdot P_{0})}{\vert\xi\vert^{2}}\cos(\gamma t\vert\xi\vert)e^{-(\alpha + \beta)\vert\xi\vert^{2}t/2}\Big\Vert +  O(t^{-\frac{n}{4}-\frac{1}{2}})\\
&=& \Vert X(t,\cdot) + Y(t,\cdot)\Vert + O(t^{-\frac{n}{4}-\frac{1}{2}})\nonumber\\
&\geq& \Vert X(t,\cdot)\Vert -\Vert Y(t,\cdot)\Vert + O(t^{-\frac{n}{4}-\frac{1}{2}}),\nonumber
\end{eqnarray}
where we have just defined as
\[X(t,\xi) := (i\xi)\frac{\sin(\gamma t\vert\xi\vert)}{\vert\xi\vert}Q_{0}e^{-(\alpha + \beta)\vert\xi\vert^{2}t/2},\]
\[
Y(t,\xi) := - \{P_{0} - \frac{\xi(\xi\cdot P_{0})}{\vert\xi\vert^{2}}\}e^{-\alpha\vert\xi\vert^{2}t} - \frac{\xi(\xi\cdot P_{0})}{\vert\xi\vert^{2}}\cos(\gamma t\vert\xi\vert)e^{-(\alpha + \beta)\vert\xi\vert^{2}t/2}.
\]

Now we will prove the following lemma.
\begin{lem}\label{lem1.3} Let $n \geq 1$.  Then, there exist three real numbers $C_{j} > 0$ {\rm (}$j= 1,2,3${\rm )}  depending on $n$, $\alpha$ or $\alpha+\beta$, such that\\
{\rm (1)}\,$C_{1}^{2}\vert P_{0}\vert^{2}t^{-\frac{n}{2}} \leq\displaystyle{\int_{{\bf R}_{\xi}^{n}}}\Big\vert\displaystyle{\frac{\xi(\xi\cdot P_{0})}{\vert\xi\vert^{2}}}e^{-\alpha\vert\xi\vert^{2}t}\Big\vert^{2}d\xi \leq C_{1}^{-2}\vert P_{0}\vert^{2}t^{-\frac{n}{2}}$,\\
{\rm (2)}\,$C^{2}_{2}t^{-\frac{n}{2}} \leq \displaystyle{\int_{{\bf R}_{\xi}^{n}}}\Big\vert (i\xi)e^{-\frac{(\alpha + \beta)\vert\xi\vert^{2}t}{2}}\displaystyle{ \frac{\sin(\gamma t\vert\xi\vert)}{\vert\xi\vert}}\Big\vert^{2}d\xi \leq C_{2}^{{-2}}t^{-\frac{n}{2}}$,\\
{\rm (3)}\,$C_{3}^{2}\vert P_{0}\vert^{2}t^{-\frac{n}{2}} \leq \displaystyle{\int_{{\bf R}_{\xi}^{n}}}\Big\vert\displaystyle{\frac{\xi(\xi\cdot P_{0})}{\vert\xi\vert^{2}}}e^{-\frac{(\alpha + \beta)\vert\xi\vert^{2}t}{2}}\cos(\gamma t\vert\xi\vert)\Big\vert^{2}d\xi \leq C_{3}^{-2}\vert P_{0}\vert^{2}t^{-\frac{n}{2}}$,\\
for large $t \gg 1$.
\end{lem} 

\begin{rem}\label{rem3.1}
{\rm It is well known that $C_{4}t^{-\frac{n}{4}} \leq \Vert e^{-\alpha\vert\xi\vert^{2}t}\Vert \leq C_{4}^{-1}t^{-\frac{n}{4}}$ as $t \to +\infty$ with some constant $C_{4} >0$  depending on $\alpha>0$ and the dimension $n$.}
\end{rem}

Let us postpone the proof of Lemma 3.1 for a while. Once Lemma 3.1 could be proved, one can proceed the proof of Theorem 1.3 as follows.\\
{\it Proof of Theorem 1.3.} First, from Lemma 3.1 one can get
\begin{equation}
\Vert X(t,\cdot)\Vert \geq C_{2}\vert Q_{0}\vert t^{-\frac{n}{4}},
\end{equation}
\begin{equation}
\Vert Y(t,\cdot)\Vert \leq C_{*}\vert P_{0}\vert t^{-\frac{n}{4}},
\end{equation}
with some constant $C_{*} > 0$ which depends on  $C_{j} > 0$ ($j =  1,3,4$). Thus, it follows from (3.1), (3.2) and (3.3) that 
\[\Vert v(t,\cdot)\Vert \geq \Big(C_{2}\vert Q_{0}\vert-C_{*}\vert P_{0}\vert\Big)t^{-\frac{n}{4}} + o(t^{-\frac{n}{4}}).\]
So, if $\vert Q_{0}\vert \ne 0$, and $\displaystyle{\frac{C_{2}}{2}}\vert Q_{0}\vert > (\displaystyle{\frac{C_{2}}{2}} + C_{*})\vert P_{0}\vert$, then as $t \to +\infty$ one can get
\[\Vert v(t,\cdot)\Vert \geq \frac{1}{2}\Big(C_{2}\vert Q_{0}\vert-C_{*}\vert P_{0}\vert\Big)t^{-\frac{n}{4}} \geq \frac{C_{2}}{4}(\vert Q_{0}\vert + \vert P_{0}\vert)t^{-\frac{n}{4}},\]
which implies the desired estimate from below. The estimate from above is a direct consequence of Lemma 3.1, Remark 3.1, and Theorem 1.2.
\par
\hfill
$\Box$
\par
 
\noindent
{\it Proof of Lemma 3.1.}\, The estimate from above of item (1) easily follows from Remark \ref{rem3.1} (see also (\ref{elem-est}) ). The estimate from above of  item (3) is estimated as follows because of the Schwarz inequality,
\[\displaystyle{\int_{{\bf R}_{\xi}^{n}}}\vert\displaystyle{\frac{\xi(\xi\cdot P_{0})}{\vert\xi\vert^{2}}}e^{-(\alpha + \beta)\vert\xi\vert^{2}t/2}\cos(\gamma t\vert\xi\vert)\vert^{2}d\xi \]
\[\leq \vert P_{0}\vert^{2}\displaystyle{\int_{{\bf R}_{\xi}^{n}}} e^{-(\alpha + \beta)\vert\xi\vert^{2}t}d\xi \leq C_{3}^{-2}\vert P_{0}\vert^{2}t^{-\frac{n}{2}}.\]
Concerning the estimate from above of (2) is also an easy exercise, so we omit its check (see Remark 3.1).

About the estimate from below of (2), one follows a device from \cite{Ik-0}. For this check, we set
\[I(t) := \displaystyle{\int_{{\bf R}_{\xi}^{n}}}\vert (i\xi)e^{-(\alpha + \beta)\vert\xi\vert^{2}t/2}\displaystyle{ \frac{\sin(\gamma t\vert\xi\vert)}{\vert\xi\vert}}\vert^{2}d\xi.\]  
Then, easily one can get a series of equalities below with the help of polar co-ordinate transform
\[I(t) = (\int_{\vert\omega\vert = 1}d\omega)(\alpha + \beta)^{-\frac{n}{2}}t^{-\frac{n}{2}}\int_{0}^{\infty}e^{-\theta^{2}}\theta^{n-1}\sin^{2}(\frac{\gamma}{\sqrt{\alpha + \beta}}\sqrt{t}\theta)d\theta\]
\[=  \frac{S_{0}}{2}(\int_{\vert\omega\vert=1}d\omega)(\alpha + \beta)^{-\frac{n}{2}}t^{-\frac{n}{2}}\]
\begin{equation}
-\frac{1}{2}(\int_{\vert\omega\vert=1}d\omega)(\alpha + \beta)^{-\frac{n}{2}}t^{-\frac{n}{2}}\int_{0}^{\infty}e^{-\theta^{2}}\theta^{n-1}\cos(\frac{2\gamma}{\sqrt{\alpha + \beta}}\sqrt{t}\theta)d\theta,
\end{equation}
where
\[S_{0} = \int_{0}^{\infty}e^{-\theta^{2}}\theta^{n-1}d\theta.\]
Since
\[\lim_{t \to +\infty}\int_{0}^{\infty}e^{-\theta^{2}}\theta^{n-1}\cos(\frac{2\gamma}{\sqrt{\alpha + \beta}}\sqrt{t}\theta)d\theta = 0,\]
because of the fact $e^{-\theta^{2}}\theta^{n-1} \in L^{1}(0,\infty)$ and the Riemann-Lebesgue Theorem, it follows from (3.4) that 
\begin{equation}
I(t) = \frac{S_{0}}{2}(\int_{\vert\omega\vert=1}d\omega)(\alpha + \beta)^{-\frac{n}{2}}t^{-\frac{n}{2}} + o(t^{-\frac{n}{2}})
\end{equation}
as $t \to +\infty$. Thus, one can get 
\[I(t) \geq \frac{S_{0}}{4}(\int_{\vert\omega\vert=1}d\omega)(\alpha + \beta)^{-\frac{n}{2}}t^{-\frac{n}{2}},\quad t \gg 1,\]
which implies the desired estimate from below. 

The estimate from below of item (3) comes from the same idea used in \cite[(2.4)]{IO}. In fact, if one sets the conical region $K \subset {\bf R}_{\xi}^{n}$ as 

\[
K:=\left\{\xi\in{\bf R}_{\xi}^n \ \bigg\vert\ \frac{\xi}{\vert\xi\vert}\cdot\frac{P_{0}}{\vert P_{0}\vert}\geq \frac12\right\},
\]
then one can observe that
\[\int_{{\bf R}^n}\vert\frac{\xi}{\vert\xi\vert}\cdot P_{0}\vert^{2}e^{-t(\alpha+\beta)\vert\xi\vert^{2}}\left\vert\cos(\gamma t\vert\xi\vert)\right\vert^2\,d\xi
\geq
\frac{\vert P_{0}\vert^2}{4}\int_{K}e^{-t(\alpha+\beta)\vert\xi\vert^{2}}\left\vert\cos(\gamma t\vert\xi\vert)\right\vert^2\,d\xi\]
\[= \frac{\vert P_{0}\vert^2}{4}\int_0^\infty\int_{K\cap\{\vert\xi\vert=r\}}e^{-t(\alpha+\beta)r^{2}}\left\vert\cos(\gamma tr)\right\vert^2\,dS\,dr\]
\[= \frac{c(n)\vert P_{0}\vert^2}{4}\int_0^\infty r^{n-1}e^{-t(\alpha+\beta)r^{2}}\left\vert\cos(\gamma tr)\right\vert^2\,dr\] 
\[= \frac{c(n)\vert P_{0}\vert^2}{8}\left(\int_0^\infty r^{n-1}e^{-t(\alpha+\beta)r^{2}}\,dr + \int_0^\infty r^{n-1}e^{-t(\alpha+\beta)r^{2}}\cos(2\gamma tr)\,dr\right),\] 
where the constant $c(n)>0$ is the area of $K\cap\{\vert\xi\vert=1\}$. So, one can proceed similar computations to item (2) in order to check the desired estimate from below of (3) by using the Riemann-Lebesgue Theorem. Finally, the estimate from below of item (1) is derived by the same argument as just above more easily, and for this we shall omit its detail. 
\par
\hfill
$\Box$
\par

\vspace{0.1cm}
\noindent{\em Acknowledgment.}
\smallskip
The work of the second author was supported in part by Grant-in-Aid for Scientific Research (C)15K04958 of JSPS. 



\begin{thebibliography}{99}
\bibitem{BK} J. Brezina and Y. Kagei, Decay properties of solution to the linearized compressible Navier-Stokes equation around time-periodic parallel flow, Math. Models Methods Appl. Sci. 22 (2012), no. 7, 1250007, 53 pp.

\bibitem{CR} S. Chowdhury and M. Ramaswamy, Optimal control of linearized compressible Navier-Stokes equations, ESAIM Control Optim. Calc. Var. 19, no. 2 (2013), 587-615.

\bibitem{RCR2} R. Coimbra Char\~ao, C. R. da Luz, R. Ikehata, New decay rates for a problem of plate dynamics with fractional damping, J. Hyperbolic Differ. Equ. 10, no. 3 (2013), 563-575.

\bibitem{D} K. Deckelnick, $L^{2}$-decay for the compressible Navier-Stokes equations in unbounded domains, Comm. Partial Diff. Eqns 18 (1993), 1445-1476.

\bibitem{HZ} D. Hoff and K. Zumbrun, Multi-dimensional diffusion waves for the Navier-Stokes equations of compressible flow, Indiana Univ. Math. J. 44, no.2 (1995), 603-676.
\bibitem{HZ-1} D. Hoff and K. Zumbrun, Pointwise decay estimates for multidimensional Navier-Stokes diffusion waves, Z. Angew. Math. Phys. 48 (1997), 1-18.

\bibitem{Ik-3} R. Ikehata, New decay estimates for linear damped wave equations and its application to nonlinear problem, Math. Meth. Appl. Sci. 27 (2004), 865-889.

\bibitem{Ik-0} R. Ikehata, Asymptotic profiles for wave equations with strong damping, J. Diff. Eqns 257 (2014), 2159-2177.

\bibitem{IKM} R. Ikehata, T. Kobayashi and T. Matsuyama, Remark on the $L_{2}$ estimates of the density for the compressible Navier-Stokes flow in ${\bf R}^{3}$, Nonlinear Analysis, TMA 47 (2001), 2519-2526.

\bibitem{IO} R. Ikehata and M. Onodera, Remarks on large time behavior of the $L^{2}$-norm of solutions to strongly damped wave equations, submitted (2016).

\bibitem{ITY} R. Ikehata, G. Todorova and B. Yordanov, Wave equations with strong damping in Hilbert spaces, J. Diff. Eqns 254 (2013), 3352-3368.

\bibitem{Kawa} S. Kawashima, System of a hyperbolic-parabolic composite type, with applications to the equations of magnetohydrodynamics, Ph.D. thesis, Kyoto University (1983).

\bibitem{KS} T. Kobayashi and Y. Shibata, Remark on the rate of decay of solutions to linearized compressible Navier-Stokes equations, Pacific J. Math. 207, no.1 (2002), 199-234.

\bibitem{LN} Tai-P. Liu and S. E. Noh, Wave propagation for the compressible Navier-Stokes equations, J. Hyperbolic Differ. Equ. 12, no. 2 (2015), 385-445.
\bibitem{TW} Tai-P. Liu and W. Wang, The pointwise estimates of diffusion wave for the Navier-Stokes systems in odd multi-dimensions, Comm. Math. Phis. 196 (1998), 145-173.

\bibitem{LIC} C. R. da Luz, R. Ikehata and R. C. Char\~ao, Asymptotic behavior for abstract evolution differential equations of second order, J. Diff. Eqns 259 (2015), 5017-5039.

\bibitem{MW} S. Ma and J. Wang, Decay rates to viscous contact waves for the compressible Navier-Stokes equations, Journal of Math. Phys. 57, 021501 (2016).

\bibitem{MN} A. Matsumura and T. Nishida, The initial value problem for the equations of motion of compressible viscous and heat-conductive gases, J. Math. Kyoto Univ. 20-1 (1980), 67-104.


\bibitem{p} G. Ponce, Global existence of small solutions to a class of nonlinear evolution equations, Nonlinear Anal. 9(5) (19), 399-418.

\bibitem{shibata} Y. Shibata, On the rate of decay of solutions to linear viscoelastic equation, Math. Meth. Appl. Sci. 23 (3) (2000), 203-226.

\end{thebibliography}
\end{document}